\documentclass[11pt]{amsart}
\usepackage{amssymb}
\newcommand{\bbC}{{\mathbb C}}
\newcommand{\bbZ}{{\mathbb Z}}
\newcommand{\bbR}{{\mathbb R}}
\newcommand{\bbQ}{{\mathbb Q}}
\newcommand{\calV}{{\mathcal V}}

\DeclareMathOperator{\Spec}{Spec}
\DeclareMathOperator{\HH}{H}
\DeclareMathOperator{\RR}{R}
\DeclareMathOperator{\CH}{CH}

\DeclareMathOperator{\MHM}{MHM}

\DeclareMathOperator{\tr}{tr}
\DeclareMathOperator{\Sym}{Sym}
\DeclareMathOperator{\IJ}{IJ}
\DeclareMathOperator{\Ext}{Ext}
\DeclareMathOperator{\Div}{div}
\DeclareMathOperator{\gr}{gr}
\DeclareMathOperator{\cl}{cl}
\DeclareMathOperator{\dR}{dR}
\DeclareMathOperator{\Db}{Db}
\DeclareMathOperator{\ADb}{ADb}
\DeclareMathOperator{\image}{im}

\newcommand{\into}{\hookrightarrow}
\newcommand{\by}[1]{\xrightarrow{#1}}
\newcommand{\ext}[1]{\operatorname{\stackrel{#1}{\wedge}}}
\newcommand{\tensor}{\otimes}
\newcommand{\isom}{\cong}
\newcommand{\abuts}{\Rightarrow}
\theoremstyle{plain}
\newtheorem{lemma}{Lemma}
\newtheorem{conj}{Conjecture}
\newtheorem{sublem}{Sublemma}
\theoremstyle{definition}

\newtheorem{defn}{Definition}
\newtheorem{problem}{Problem}
\newenvironment{diagram}[1]{\arraycolsep=\doublerulesep\begin{array}{#1}
    }{\end{array}}

\begin{document}

\title{Higher Abel-Jacobi Maps}

\author[J. Biswas]{Jishnu Biswas}
\address{IMSc, CIT Campus, Tharamani, Chennai 600 113, India.}
\email{jishnu@imsc.ernet.in}

\author[G. Dayal]{Gautham Dayal}
\address{CMI, G N Chetty Road, T Nagar, Chennai 600 017, India.}
\email{gautham@imsc.ernet.in}

\author[K. H. Paranjape]{Kapil H.~Paranjape}
\address{IMSc, CIT Campus, Tharamani, Chennai 600 113, India.}
\email{kapil@imsc.ernet.in}

\author{G.~V.~Ravindra}
\address{IMSc, CIT Campus, Tharamani, Chennai 600 113, India.}
\email{ravindra@imsc.ernet.in}

\maketitle

\section*{Introduction}

We work over a subfield $k$ of $\bbC$, the field of complex
numbers.  For a smooth variety $V$ over $k$, the Chow group of
cycles of codimension $p$ is defined (see \cite{Fulton}) as
\[ \CH^{p}(V)=\frac{Z^{p}(V)}{R^{p}(V)}\]
where the group of cycles $Z^{p}(V)$ is the free abelian group on
scheme-theoretic points of $V$ of codimension $p$ and rational
equivalence $R^{p}(V)$ is the subgroup generated by cycles of the form
$\Div_{W}(f)$ where $W$ is a subvariety of $V$ of codimension $(p-1)$
and $f$ is a non-zero rational function on it. There is a natural
cycle class map
\[ \cl_{p}:\CH^{p}(V)\to \HH^{2p}(V)\]
where the latter denotes the singular cohomology group
$\HH^{2p}(V(\bbC),\bbZ)$ with the (mixed) Hodge structure given by
Deligne (see \cite{HodgeII}). The kernel of $\cl_{p}$ is denoted by
$F^{1}\CH^{p}(V)$. There is an Abel-Jacobi map (see
\cite{Griffiths}),
\[ \Phi_{p}:F^{1}\CH^{p}(V)\to \IJ^{p}(\HH^{2p-1}(V))\]
where the latter is the intermediate Jacobian of a Hodge
structure, defined as follows
\[ \IJ^{p}(H)=\frac{H\tensor \bbC}{F^{p}(H\tensor \bbC)+H}.\]
We note for future reference that for a pure Hodge structure of weight
$2p-1$ (such as the cohomology of a smooth \emph{projective}
variety) we have the natural isomorphism,
\[ \frac{H\tensor \bbR}{H}\widetilde{\to}\IJ^{p}(H)\]
The kernel of $\Phi_{p}$ is denoted by $F^{2}\CH^{p}(V)$.

The conjecture of S.~Bloch says (see \cite{Jannsen}) that there is a
filtration $F^{\cdot}$ on $\CH^{p}(V)$ which extends the $F^{1}$ and $F^{2}$
defined above.  Moreover, the associated graded group
$\gr^{k}_{F}\CH^{p}(V)$ is \emph{governed} by the cohomology 
groups
$\HH^{2p-k}(V)$ for each integer $k$ (upto torsion). More precisely,
if $N^{l}\HH^{m}(V)$ denotes the filtration by co-niveau (see
\cite{Grothendieck}) which is generated by cohomology classes
supported on subvarieties of codimension $\leq l$, then
$\gr^{k}_{F}\CH^{p}(V)$ is actually governed by the quotient 
groups, \[ \HH^{2p-k}(V)/N^{p-k+1}\HH^{2p-k}(V).\]
Specifically, in the case when $V$ is a smooth projective surface with
geometric genus 0 (so that $\HH^{2}(V)=N^{1}\HH^{2}(V)$) this
conjecture implies that $F^{2}\CH^{2}(V)$ is torsion (and thus 0 by a
theorem of Roitman \cite{Roitman}).

The traditional Hodge-theoretic approach to study this problem is
based on the fact that the intermediate Jacobian $\tensor\bbQ$ can be
interpreted as the extension group $\Ext^{1}(\bbQ(-p),H)$ in the
category of Hodge structures. One can then propose that the associated
graded groups $\gr_F^k\CH^p(V)\tensor \bbQ$ should be interpreted as
the higher extension groups $\Ext^{k}(\bbQ(-p),\HH^{2p-k}(V))$ for
$k\geq 2$.  Unfortunately, there are no such extension groups in the
category of Hodge structures. Thus it was proposed that 
all these extension groups be computed in a suitable category of
mixed motives.\footnote{%
  Such a category has recently been constructed by M.~V.~Nori
  (unpublished).}

Even if such a category is constructed a Hodge-theoretic
interpretation of these extension groups would be useful. In section~2
we discuss M.~Green's approach (see \cite{Green}) called the Higher
Abel-Jacobi map. In section~3 we provide a counter-example to show
that Green's approach does not work; a somewhat more complicated
example was earlier obtained by C.~Voisin (see \cite{Voisin}).  In
section~4 we introduce an alternative approach based on
Deligne-Beilinson cohomology and its interpretation in terms of
Morihiko Saito's theory of Hodge modules; such an approach has also
been suggested earlier by M. Asakura and independently by M. Saito
(see \cite{Asakura} and \cite{Saito}). Following this approach it
becomes possible to deduce Bloch's conjecture from some conjectures of
Bloch and Beilinson for cycles and varieties defined over a number
field (see \cite{Jannsen}).

\section{Green's Higher Abel-Jacobi Map}

The fundamental idea behind M.~Green's construction can be interpreted
as follows (see \cite{Voisin}). One expects that the extension groups
are \emph{effaceable} in the abelian category of mixed motives. Thus
the elements of $\Ext^{k}$ can be written in terms of $k$ different
elements in various $\Ext^{1}$'s.  The latter groups can be understood
in terms of Hodge theory, via the Intermediate Jacobians. So we can
try to write the $\Ext^k$ as a sub of a quotient of a (sum of) tensor
products of Intermediate Jacobians.

Specifically, consider the case of a surface $S$. Let $C$ be a curve,
then we have a product map (see \cite{Fulton}), 
\[ \CH^{1}(C)\times \CH^{2}(C\times S)\to \CH^{2}(S)\]
which in fact respects the filtration $F^{\cdot}$ (see
\cite{Jannsen}), so that we have
\[ F^{1}\CH^{1}(C)\times F^{1}\CH^{2}(C\times S)\to F^{2}\CH^{2}(S).\]
Conversely, we can use an argument of Murre (see \cite{Murre}) to
show,
\begin{lemma}
  Given any cycle class $\xi$ in $F^{2}\CH^{2}(S)$ there is a
  curve $C$ so that $\xi$ is in the image of the map,
\[ F^{1}\CH^{1}(C)\times F^{1}\CH^{2}(C\times S)\to F^{2}\CH^{2}(S).\]
\end{lemma}
\begin{proof}
  Let $z$ be a cycle representing the class $\xi$. There is a smooth
  (see \cite{CPR}) curve $C$ on $S$ that contains the support of $z$.
  Hence it is enough to show that there is a homologically trivial
  cycle $Y$ on $C\times S$ so that $(z,Y)\mapsto z$ for every cycle
  $z$ on $C$ such that the image under $\CH^1(C)\to\CH^2(S)$ lies in
  $F^2\CH^2(S)$. Let $\Gamma$ denote the graph of the inclusion
  $\iota:C\into S$. Then clearly $(z,\Gamma)\mapsto z$ but $\Gamma$ is not
  homologically trivial.
  
  Choose a point $p$ on $C$. Now, by a result of Murre (see
  \cite{Murre}), for some positive integer $m$ we have an expresssion
  in $\CH^2(S\times S)$
  \[
   m \Delta_S = m(p\times S + S \times p) + X_{2,2} + X_{1,3} + X_{3,1}
  \] 
  where $\Delta_S$ is the diagonal and $X_{i,j}$ is a cycle so that
  its cohomology class has non-zero K\"unneth component only in
  $\HH^i(S)\tensor\HH^j(S)$. In particular, $X_{1,3}$ gives a map
  $F^1\CH^2(S)\to F^1\CH^2(S)$ which induces multiplication by $m$ on
  $\IJ^2(\HH^3(S))$. Since $p\times S$ and $S\times p$ induce $0$ on
  $F^1\CH^2(S)$ it follows that the correspondence $X_{2,2}+X_{3,1}$
  induces multiplication by $m$ on $F^2\CH^2(S)$.
  
  Now, $\Gamma=(\iota\times 1_S)^*(\Delta_S)$ so we have an expression
  \[ m \Gamma = m C \times p + (\iota\times 1_S)^* X_{2,2} +
            (\iota\times 1_S)^* X_{1,3}
  \]
  Let $D=p_{2*}((\iota\times 1_S)^* X_{2,2})$. Then the cohomology
  class of $Y=(\iota\times 1_S)^* X_{2,2} - p\times D$ is 0. Moreover,
  the map $F^1\CH^1(C)\to F^1\CH^2(S)$ induced by $p\times D$ is zero.
  Thus, by the above propery of $X_{2,2}+X_{3,1}$ we see that
  $mz=(z,m\Gamma)=(z,Y)$ for any $z$ in $F^1\CH^1(C)$ whose image lies
  in $F^2\CH^2(S)\tensor\bbQ$. By Roitman's theorem (see
  \cite{Roitman}) the group $F^2\CH^2(S)$ is divisible.  Hence, we
  conclude the result.
\end{proof}
We now use the Abel-Jacobi maps to interpret the two terms on the
left-hand side in terms of Hodge theory.

Firstly, we have the classical Abel-Jacobi isomorphisms
$F^{1}\CH^{1}(C)=J(C)=\IJ^{1}(\HH^{1}(C))$. Let
$\HH^{2}(S)_{\tr}=\HH^{2}(S)/N^{1}\HH^{2}(S)$ denote the lattice of
transcendental cycles on $S$. Consider the factor
$\IJ^{2}(\HH^{1}(C)\tensor\HH^{2}(S)_{\tr})$ of the intermediate
Jacobian $\IJ^2(\HH^3(C\times S))$. We can compose the Abel-Jacobi map
with the projection to this factor to obtain 
\[ 
 F^{1}\CH^{2}(C\times S)\to
 \IJ^{2}(\HH^{1}(C)\tensor\HH^{2}(S)_{\tr}).
\]
Using the identification $\IJ^{p}(H)=H\tensor (\bbR /\bbZ)$ for a pure
Hodge structure $H$ of weight $2p-1$ we have
\[
   \IJ^{1}(\HH^{1}(C))\tensor \IJ^{2}(\HH^{1}(C)\tensor
                      \HH^{2}(S)_{\tr})=
    \HH^{1}(C)^{\tensor 2}\tensor \HH^{2}(S)\tensor
                           (\bbR /\bbZ)^{\tensor 2}
\]
The pairing $\HH^{1}(C)^{\tensor 2}\to \HH^{2}(C)=\bbZ$,
given by the cup product, can be used to further collapse the latter
term. Thus, we obtain a diagram,
\[
 \begin{diagram}{ccc}
       F^{1}\CH^{1}(C)\times F^{1}\CH^{2}(C\times S) &
                    \to  & F^{2}\CH^{2}(S)\\
          \downarrow  &  & \\
       \IJ^1(\HH^{1}(C))\times \IJ^2(\HH^{1}(C)\tensor
          \HH^{2}(S)_{\tr}) &
                    \to  & \HH^{2}(S)_{\tr}\tensor
                                        (\bbR /\bbZ)^{\tensor 2}
 \end{diagram}
\]
\begin{defn}
  Green's second intermediate Jacobian $J^{2}_{2}(S)$ is defined
  as the universal push-out of all the above diagrams as $C$ is
  allowed to vary. The Higher Abel-Jacobi map is defined as
  the natural homomorphism 
  \[ \Psi^{2}_{2}:F^{2}\CH^{2}(S)\to J^{2}_{2}(S). \]
\end{defn}
By the above lemma it follows that $J^{2}_{2}(S)$ is a quotient of
$\HH^{2}(S)_{\tr}\tensor (\bbR /\bbZ)^{\tensor 2}$. The question is
whether this constructs the required $\Ext^2$.
\begin{problem}[Green]
  Is $\Psi^{2}_{2}$ injective?
\end{problem}

\section{Non-injectivity of Green's Map}

We now compute Green's Higher Abel-Jacobi map for the case of a
surface of the form $\Sym^2(C)$, where $C$ is a smooth projective
curve. Using this we show that this map is not injective when $C$ is a
curve of genus at least two whose Jacobian is a simple abelian
variety. 

\begin{lemma}
  Let $Z\in \CH^{2}(D \times C \times S)$ be a cycle, where $D$, $C$
  are smooth curves and $S$ a smooth surface. Then we have a
  commutative diagram
\[
\begin{diagram}{ccc}
F^{1}\CH^1(D) \tensor F^{1}\CH^1(C) &
           {\by{p_{3*}(p_{12}^*(\_)\cdot Z)}}
                      & {F^{2}\CH^{2}(S)} \\
\downarrow && \downarrow \\
\IJ^1(H^1(D)) \tensor \IJ^1(H^1(C)) &
  \by{z} ~~~ \HH^{2}(S)_{\tr}\tensor (\bbR/\bbZ)^{\tensor 2} ~~~
       \to  & J^{2}_{2}(S) 
\end{diagram}
\]
Here the map $z$ is the composite as follows. The cohomology class of
$Z$ gives a map $\HH^1(D)\tensor\HH^1(C)\to\HH^2(S)$; we further
project to $\HH^2(S)_{\tr}$. Now tensor with
$(\bbR/\bbZ)^{\tensor 2}$ and identify the resulting left-hand term
with the product of the Intermediate Jacobians.
\end{lemma}
We note that the vertical arrow on the left is an isomorphism.
\begin{proof}
  By the functoriality of the Abel-Jacobi map we have a commutative
  diagram
\[
  \begin{diagram}{ccc}
     F^{1}\CH^1(D) & \by{p_{23*}(p_1^*(\_)\cdot Z)} &
                                  F^{1}\CH^{2}(C\times S) \\
           \downarrow && \downarrow \\
     \IJ^1(\HH^1(D)) &
           \by{1_{(\bbR/\bbZ)}\tensor p_{23*}(p_1^*(\_)\cup [Z])}
         & \IJ^2(\HH^{1}(C)\tensor \HH^{2}(S))
  \end{diagram}
\]
By projection we can replace the bottom right corner with
$\IJ^2(\HH^{1}(C)\tensor\HH^{2}(S)_{\tr})$. Now we tensor this with
the Abel-Jacobi map for $C$ to obtain,
\[
\begin{diagram}{ccrcl}
F^{1}\CH^1(D)\tensor F^{1}\CH^1(C) & \to &
            F^1\CH^1(C) & \tensor &F^{1}\CH^{2}(C\times S) \\
\downarrow && & \downarrow \\
     \IJ^1(\HH^1(D)) \tensor \IJ^1(\HH^1(C))  & \to &
                \IJ^1(\HH^1(C)) & \tensor & \IJ^2(\HH^{1}(C)\tensor \HH^{2}(S)_{\tr})
\end{diagram}
\]
The required commutative diagram now follows from the definition of
$J^2_2(S)$. 
\end{proof}
We now apply this lemma to the case $C=D$ and $S=\Sym^{2}(C)$. In this
case we take $Z$ to be the graph of the quotient morphism $q:C\times
C\to \Sym^2(C)$. We then compute that the cohomological
correspondence given by $[Z]$ factors as
\[
  \HH^1(C)\tensor\HH^1(C) \to \ext{2}\HH^{1}(C) \to
                  \HH^{2}(\Sym^2(C))_{\tr}
\]
By the above lemma we obtain a factoring,
\[
\begin{diagram}{ccc}
F^{1}\CH^1(C) \tensor F^{1}\CH^1(C) &
           \by{p_{3*}(p_{12}^*(\_)\cdot Z)}
                      & F^{2}\CH^{2}(\Sym^2(C)) \\
\downarrow && \downarrow \\
(H^1(C) \tensor \bbR/\bbZ)^{\tensor 2}
      & \to ~~~ \ext{2}\HH^{1}(C)\tensor (\bbR/\bbZ)^{\tensor 2}
       ~~~ \to & J^{2}_{2}(\Sym^2(C)) 
\end{diagram}
\]
The image of the tensor product of a pair of elements of
$\IJ^1(\HH^1(C))$ of the form $v\tensor \alpha$ and $v\tensor \beta$
must therefore be 0 in $J^{2}_{2}(\Sym^2(C))$.

The description of $F^2\CH^2(\Sym^2(C))$ is given by the following
lemma that is similar to one in \cite{Bloch},
\begin{lemma}
  The homomorphism 
  \[
    Z_{*}: F^{1}\CH^1(C) \tensor F^{1}\CH^1(C)
                     \to F^{2}\CH^{2}(\Sym^2(C))
  \]
  is surjective.
\end{lemma}
\begin{proof}
The following composite  map  is multiplication by  2
\[  F^2\CH^2(\Sym^2(C)) \by{q^*} F^2\CH^2(C\times C)
                   \by{q_*}   F^2\CH^2(\Sym^2(C))
\]
By the divisibility of $F^2\CH^2(S)$  for  a surface $S$  we see that
the lemma follows from the following result.
\end{proof}
\begin{sublem}
  Fix a base point $p$ on $C$. Then the filtration $F$ of
  $\CH^2(C\times C)$ is explicitly described as follows
\begin{multline*}
  F^2\CH^2(C\times C) = \image(J(C)\tensor J(C)) \subset \\
     F^1\CH^2(C\times C)  =   F^2\CH^2(C\times C) + \image(J(C)\times p)
                          +  \image(p\times J(C))  \\
       \subset \CH^2(C\times C)  =  F^1\CH^2(C\times C)  +  \bbZ\cdot (p,p)
\end{multline*}
\end{sublem}
\begin{proof}
  Let $a$, $b$ be points on $C$; we get points $[a-p]$ and $[b-p]$
  of $J(C)$. The image of $[a-p]\tensor[b-p]$ in $\CH^2(C\times C)$ is
  $(a,b)+(p,p)-(a,p)-(p,b)$. Thus, we have an   expression
  \[  (a,b)  =   \image([a-p]\tensor[b-p]) +
                     \image([a-p]\times p) +
                     \image(p\times [b-p]) +
                           (p,p)
  \]
  Now, any cycle $\xi$ in $F^1\CH^2(C\times C)$ can be written as
  $\sum_{i=1}^n (a_i,b_i) - n\cdot (p,p)$. The Albanese variety of
  $C\times C$ is $J(C)\oplus J(C)$ and the image of $\xi$ under the
  Albanese map is $(\sum_{i=1}^n [a_i-p],\sum_{i=1}^n [b_i-p])$. Thus,
  if the cycle is in $F^2\CH^2(C\times C)$, then $\sum_{i=1}^n [a_i-p]
  =0 = \sum_{i=1}^n [b_i-p]$. Now we combine this with the above
  expression to obtain
  \[ \xi = \sum_{i=1}^n \image([a_i-p]\tensor[b_i-p]) \]
  which proves the result.
\end{proof}
\begin{lemma}
  If $C$ is a curve of genus at least 2 such that its Jacobian variety
  is a simple abelian variety then $\Psi^{2}_{2}$ has a non-trivial
  kernel.
\end{lemma}
\begin{proof}
  By Mumford's result there are non-trivial classes in
  $F^{2}\CH^{2}(S)$.  The Jacobian variety
  $J(C)=\IJ^1(\HH^1(C))=\HH^{1}(C)\tensor \bbR /\bbZ $ is spanned by
  decomposable elements. Moreover $F^{1}\CH^1(C)\isom J(C)$. Thus
  there is a pair of elements of $F^1\CH^1(C)$ of the form $v\tensor
  \alpha $, $w\tensor \beta $ such that the image of their tensor
  product in $F^{2}\CH^{2}(S)$ is non-zero. By a result of
  Roitman, for any fixed class $f$ in $F^1\CH^1(C)$, 
  the collection 
  \[  K_f =
        \{ e\in J(C) | e \tensor f \mapsto 0 \text{~in~}
        F^{2}\CH^{2}(S) \}
  \]
  forms a countable union of abelian subvarieties of $J(C)$. Since
  $w\tensor\beta$ does not lie in $K_{v\tensor\alpha}$, the latter is
  a proper subgroup of $J(C)$. Since $J(C)$ is assumed to be simple
  this is forced to be a countable set. In particular, there is an
  element of the form $v\tensor\gamma$ which is {\em not} in
  $K_{v\tensor\alpha}$; so that the product of this
  with $v\tensor \alpha$ is non-zero in $F^{2}\CH^{2}(S)$. But we
  just saw that all such elements are mapped to 0 in $J^{2}_{2}(S)$.
\end{proof}

\section{Absolute Deligne-Beilinson Cohomology}
The fundamental idea underlying the following constructions and
definitions is as follows. A variety $V$ over $\bbC$ can be thought of
as a family of varieties over the algebraic closure
$\overline{\bbQ}\subset\bbC$ of the field of rational numbers. Even
when the variety is defined over $\bbQ$ the Chow group of such a
variety (when considered over $\bbC$) may contain cycles that are
defined over larger fields. In particular, the usual examples of
non-trivial elements in $F^2\CH^2(S)$ are defined over fields of
transcendence degree 2 (see \cite{Srinivas}). Thus, in order to detect
such cycles we must use the full force of such a ``family''-like
structure.

For any variety $V$ over $\bbC$ we consider the collection of
Cartesian diagrams
\[
 \begin{array}{ccc}
   V & \to & \calV \\
 \downarrow && \downarrow\\
 \Spec\bbC & \to & S
 \end{array}
\]
where $S$ and $\calV$ are varieties defined over $\overline{\bbQ}$,
and the lower horizontal arrow factors through the generic point of
$S$. Assume for the moment that $V$ is smooth projective, and that $S$
and $\calV$ are smooth and $\calV\to S$ is proper and smooth. Then the
relative de~Rham cohomology groups $\HH^i_{\dR}(\calV/S)$ carry the
Gauss-Manin connection; moreover, after base change to $S\tensor\bbC$
the associated local system is a variation of Hodge structure. This
has been generalised by M.~Saito (see \cite{MHM}) for all $V$ and
all choices of $S$ and $\calV$ as follows.\footnote{%
  Since we have chosen an embedding $\overline{\bbQ}\subset\bbC$ we can
  think of Hodge modules as being associated with varieties over
  $\overline{\bbQ}$ rather than with varieties over $\bbC$}
There is a
(mixed) Hodge module $\RR^i_{\dR}(\calV/S)$ on $S$ in the above
context so that its pull-back via $\Spec\bbC\to S$ is the (mixed)
Hodge structure on the cohomology of $V$.
The category $\MHM(S)$ of Hodge modules over $S$ is an abelian
category which {\em has} non-trivial $\Ext^2$'s when $S$ has dimension
at least 1. Moreover, we have a natural spectral sequence
\[
   E_1^{a,b} = \Ext^b_{\MHM(S)}(\bbQ(c),\RR^{a}_{\dR}(\calV/S))
                  \abuts \Ext^{a+b}_{\MHM(\calV)}(\bbQ(c),\bbQ)
\]
We are interested in the case $a=2p-k$, $b=k$ and $c=-p$. In this case
the latter term can be identified with the Deligne-Beilinson
cohomology $\HH^{2p}_{\Db}(\calV,\bbQ(p))$ (see \cite{Saito} and
  \cite{Beilinson}).
\begin{defn}
Let us define the absolute Deligne-Beilinson cohomology of $V$ as the
direct limit
\[           \HH^n_{\ADb}(V,\bbQ(c)) =
                  \lim_{\rightarrow} \Ext^n_{\MHM(\calV)}(\bbQ(c),\bbQ)
\]
where the limit is taken over all diagrams such as the one above.
\end{defn}
Since any algebraic cycle on $V$ (and $V$ itself) is defined over some
finitely generated field, we have
\[ \CH^p(V) = \lim_{\rightarrow} CH^p(\calV) \]
The cycle class map in Deligne-Beilinson cohomology then gives us a
cycle class map
\[ \cl^p_{\ADb}: \CH^p(V) \to \HH^{2p}_{\ADb}(V,\bbQ(p)) \]
The filtration on the latter group induced by the above spectral
sequence induces a filtration on $\CH^p(V)$. We can then ask whether
this is the filtration as required by Bloch's conjecture.

It is well known (see \cite{Schneideretal}) that the cycle class map
for Deligne-Beilinson cohomology combines the usual cycle class map to
singular cohomology with the Abel-Jacobi map. Thus, the following
conjecture implies that $\cl^p_{\ADb}$ is injective.
\begin{conj}[Bloch-Beilinson]
 If $V$ is a variety defined over a number field then
 $F^2\CH^p(V)=0$.
\end{conj}
We (of course) offer no proof of this conjecture. However, there are
examples due to C.~Schoen and M.~V.~Nori (see \cite{Schoen}),
discovered independently by M.~Green and the third author, which show
that one cannot relax the conditions in this conjecture. A
paper~\cite{GreenKap} by M.~Green and and the third author contains
these and other examples showing that $F^2\CH^p(V)$ can be non-zero
for $V$ a variety over a field of transcendence degree at least one.

\providecommand{\bysame}{\leavevmode\hbox to3em{\hrulefill}\thinspace}

\end{document}